\documentclass[11pt]{article}

\usepackage{amsmath,amssymb}
\usepackage{palatino,mathpazo}

\newtheorem{thm}{Theorem}
\newtheorem{prop}{Proposition}
\newtheorem{defn}{Definition}

\newtheorem{cor}{Corollary}

\newcommand{\sm}{\left(\begin{smallmatrix}}
\newcommand{\esm}{\end{smallmatrix}\right)}
\newcommand{\pf}{\noindent {\bf Proof.} }
\newcommand{\qed}{\hfill {$\Box$}}

\def\Z{\mathbf{Z}}
\def\Q{\mathbf{Q}}

\def\C{\mathbf{C}}
\def\H{\mathbf{H}}
\def\G{\Gamma}
\def\a{\alpha}

\def\g{\gamma}
\def\e{\varepsilon}

\def\s{\sigma}

\def\ii{i \infty}
\def\bu{\bullet}

\def\Ii
       {\widetilde{I}^{z}_{\ii} \left(
        \begin{smallmatrix}
        \a_1,&\ldots,&\a_n\\
        f_1,&\ldots,&f_n
        \end{smallmatrix}\right)}
\def\Iia
       {\widetilde{I}^{z}_{a} \left(
        \begin{smallmatrix}
        \a_1,&\ldots,&\a_n\\
        f_1,&\ldots,&f_n
        \end{smallmatrix}\right)}
\def\Fi
       {\widetilde{F}^{z}_{\ii} \left(
        \begin{smallmatrix}
        -,&\a_2,&\ldots,&\a_n\\
        f_1,&f_2,&\ldots,&f_n
        \end{smallmatrix}\right)}
\def\Fia
       {\widetilde{F}^{z}_{a} \left(
        \begin{smallmatrix}
        -,&\a_2,&\ldots,&\a_n\\
        f_1,&f_2,&\ldots,&f_n
        \end{smallmatrix}\right)}

\def\Iichi
       {\widetilde{I}^{z}_{ i} \left(
        \begin{smallmatrix}
        \a_1,&\ldots,&\a_n\\
        f_1^{\chi_1},&\ldots,&f_n^{\chi_n}
        \end{smallmatrix}\right)}

\def\Iichibar
       {\widetilde{I}^{z}_{ i} \left(
        \begin{smallmatrix}
        \a_1,&\ldots,&\a_n\\
        f_1^{\overline{\chi_1}},&\ldots,&f_n^{\overline{\chi_n}}
        \end{smallmatrix}\right)}

\def\Iic
       {\widetilde{I}^{z}_{c} \left(
        \begin{smallmatrix}
        \a_1,&\a_2,&\ldots,&\a_n\\
        f_{1},&f_2,&\ldots,&f_n
        \end{smallmatrix}\right)}

\title{Iterated period integrals and  multiple Hecke $L$-functions}
\author{YoungJu Choie{\thanks{This work was partially supported  by Priority Research Centers Program NRF
2009-0094069. }}
 \ and Kentaro Ihara\thanks{
This work was supported by Priority Research Centers Program through the National Research Foundation of Korea (NRF) funded by the Ministry of Education, Science and Technology (Grant 2009-0094069).
 \newline
Mathematics Subject Classification; Primary 11E45, Secondary 11M32.
}
} 

\newcommand{\Abstract}
{
\begin{abstract}
In this paper we express the multiple
Hecke $L$-function in terms of  a linear combination of  
iterated period integrals associated with elliptic cusp forms,  
which is introduced by Manin around 2004. 
 This expression generalizes the classical formula of Hecke 
$L$-function obtained by the Mellin transformation of a cusp form.  
Also the expression gives a way of the analytic continuation of the multiple
Hecke $L$-function. 
\end{abstract}
}

\date{}

\begin{document}

\maketitle

\Abstract

\section{Introduction}\label{sec1}

Recently, in \cite{Man1,Man2}, Manin introduces a generalization of
the period integrals of elliptic cusp forms by means of the iterated path integrals on the complex upper-half plane and discusses several
related topics. For example, he studies the iterated version of 
Mellin transformation and its functional equation, 
the properties of non-commutative generating series 
of the iterated integrals, an analogy of the classical period theory 
and its interpretation in terms of non-abelian group cohomology, and so on.
The result which we would like to focus in his paper is
an expression of the iterated period integral
in terms of special values of a multiple Dirichlet series
in their convergent region (see \S 3.2 in \cite{Man1}).

In this paper first we define the \textit{multiple
Hecke $L$-function}, which is essentially the same as the
Dirichlet series which he introduced, and show the analytic
continuation of the function. Next, as a main result, we give the
expression of the iterated period integral in terms of a linear
combination of $L$-functions, which holds in arbitrary region, and
generalizes Manin's expression. This expression also generalizes
the classical formula \eqref{mellin} below given by Mellin
transformation. As a consequence, one can write the iterated
period integral as a $\Q$-linear combination of special values of
$L$-function.

Let $\H=\{z\in \C\ |\ {\rm Im}\ z>0\}$ be the complex
upper-half plane.
Assume that $f(z)$ is a holomorphic function
on $\H$ satisfying following two conditions: \\
(i) (Fourier expansion) $f(z)$ is periodic with period one and
has a Fourier expansion
of the form $f(z)=\sum_{m=1}^\infty c_m q^m$,
where $q=e^{2\pi i z}$ whose coefficients $\{c_m\}$ have 
at most polynomial growth in $m$ when $m\rightarrow \infty$: $c_m=O(m^M)$
for some $M>0$. \\
(ii) (Cusp conditions) For any $\g \in SL_2(\Z)$, there
exists a constant $C>0$ and an integer $k$ such that
$$
(f|_k{\g})(z):=j(\g,z)^{-k}f(\g z)=O(e^{2\pi i C z}), \quad z\rightarrow \ii,
$$
where $j(\g,z)=cz+d$ and $\g z=(az+b)/(cz+d)$ for
$\g=
\left(\begin{smallmatrix}
a&b\\c&d\end{smallmatrix}\right)
$.

For example, let $\G$ be a congruence subgroup
of $SL_2(\Z)$ containing the translation matrix 
$
\left(\begin{smallmatrix}
1&1\\0&1\end{smallmatrix}\right)
$ and $S_k(\Gamma)$ be the space of 
holomorphic $\G$-cusp forms of weight $k$. 
Then $f(z)\in S_k(\Gamma)$ satisfies both (i) and (ii).

The Mellin transformation of $f(z)$ gives following expression:
\begin{align}\label{mellin}
\int_{\ii}^0 f(z)\ z^{s-1}dz=-\Gamma(s)L(f,s)
\end{align}
where  $\Gamma(s)=\int_{0}^{+\infty} e^{-t}t^{s-1}dt$
is the Gamma function and
$L(f,s):=(-2\pi i)^{-s}\sum_{m=1}^\infty{c_m}{m^{-s}}$ for
${\mathrm{Re}\ s\gg 0}$ is
the Hecke $L$-function attached to $f$ which is normalized
by $(-2\pi i)^{-s}$ for the sake of convenience.
(We take the non-positive imaginary axis as a branch cut
in the $z$-plane, so that ${\rm arg}(z)\in [-\pi/2,3 \pi/2)$ to define
the complex power throughout this paper.)
(Following \cite{Man1,Man2}, we chose $\ii$ as the base point of
paths. The minus sign in RHS of \eqref{mellin} happens for this choice.)
The condition (i) guarantees the convergence of
the $L$-function in ${\mathrm{Re}\ s\gg 0}$ and also (ii) implies
the convergence of the integral in \eqref{mellin} in $s\in \C$.
For $f\in S_k(\G)$, the special values defined by the integral
in \eqref{mellin} at critical points $s=1,\ldots,k-1$
are called \textit{periods
of $f$} and play a fundamental role in the period theory.
Shokurov \cite{Sho} constructs some varieties, so called Kuga-Sato
varieties, and interprets these numbers as periods of
relative homology of the varieties.
See also e.g. \cite{Kn, KZ} for the period theory.
The equation \eqref{mellin} gives not only the
interpretation of the periods as the special values
of $L$-function but the way of the analytic
continuation of $L$-functions.

Manin generalizes the integral in \eqref{mellin} as follows.
For $r=1,\ldots,n$, let $f_r$ be holomorphic functions on $\H$ satisfying
the conditions (i) and (ii) above (for instance $f_r\in S_{k_r}(\G)$
for integers $k_r$) and $s_r$ be
complex variables.
For points $a,z\in \overline{\H}:=\H\cup\mathbf{P}^1(\Q)$,
we fix a path
joining $a$ to $z$ on $\H$ and approaching $a$ and $z$
vertically if $a$ or $z$ is in $\mathbf{P}^1(\Q)$. Then we
consider the iterated integral along the path:
\begin{align}\label{defI}
I_a^z\left(\begin{smallmatrix}
s_1,&\ldots,&s_n\\
f_1,&\ldots,&f_n
\end{smallmatrix}\right)
:&=
\int_{a}^{z}f_1(z_1)z_1^{s_1-1}dz_1
\int_{a}^{z_1}f_2(z_2)z_2^{s_2-1}dz_2
\cdots
\int_{a}^{z_{n-1}}f_n(z_n)z_n^{s_n-1}dz_n.
\end{align}
The integral converges again because of cusp conditions of $f_r$'s
and defines a holomorphic function
in $z\in \H$ and in $(s_1,\ldots,s_n)\in \C^n$.
For a connection between these integrals and the motif theory,
see the recent paper \cite{Ichi} of Ichikawa.
The purpose of this paper is to state a generalization of \eqref{mellin} in the case of the iterated integrals.

In \S \ref{sec2}, we introduce the multiple Hecke $L$-functions (Definition \ref{defL}) and generalize the identity \eqref{mellin}
to an iterated case
(Theorem \ref{th1}). The LHS in \eqref{mellin} will be replaced
by Manin's iterated integral.
A linear combination of the multiple Hecke $L$-functions will appear in RHS.
As a consequence, one can
write the iterated integral as a $\Q$-linear
combination of special values of $L$-function.
Conversely the special values of $L$-function also can be
expressed by the sum of iterated integrals (Corollary \ref{cor1}).
In \S \ref{sec3}, we give proofs of these theorems.
In last section, we discuss the further properties of some
functions which will be introduced in \S \ref{sec2} for the proofs.

In the last part of the Introduction, it is worth mentioning that
the integral \eqref{defI} can be seen
as an analogy of the \textit{multiple polylogarithm\ }(MPL):
\begin{align}\label{mzv}
Li_{n_1,\ldots,n_r}(z):=
\sum_{m_1>\cdots>m_r>0}
\frac{z^{m_1}}{m_1^{n_1}\cdots  m_r^{n_r}}
=\int_0^z\eta_1
\eta_2
\cdots
\eta_n,
\end{align}
where $|z|<1$, $n_l$ are positive integers with $n=\sum n_l$ and
$\eta_l(z)$
are holomorphic $1$-forms on $\mathbf{P}^1(\C)\setminus \{0,1,\infty\}$
defined by
$\eta_l(z)=dz/(1-z)$ or $dz/z$
if $l=\sum_{m=1}^{j}n_m$ for some $j=1,\ldots,r$
or otherwise, respectively.
If $n_1>1$, the limits of \eqref{mzv} as
$z\rightarrow 1$ exist and are called
\textit{multiple zeta values} (MZVs).
It is known that MZVs have a rich theory and
a geometric origin associated to
$\mathbf{P}^1(\C)\setminus \{0,1,\infty\}$
or the mixed Tate motifs over $\Z$.
See \cite{Gon1,Gon2,Z}.

As mentioned in \cite{Man1,Man2}, the iterated integral \eqref{defI} 
can be seen as an analogy of \eqref{mzv}. 
The $1$-forms on  $\H$ replaced with those on 
$\mathbf{P}^1(\C)\setminus \{0,1,\infty\}$.
In this reason, we expect a rich
theory behind the values defined by \eqref{defI}
such as MZVs. As a consequence of results in this paper,
those values will be able to interpreted as the special
values of multiple Hecke $L$-function. For a related topic, see 
\cite{I}.

\section{Periods and multiple Hecke $L$}\label{sec2}

In this section we state our main results. The proofs will be
given in $\S \ref{sec3}$.
Recall that $f_r$ ($r=1,\ldots,n$) are holomorphic functions on $\H$ satisfying the conditions (i) and (ii) in $\S \ref{sec1}$.
Suppose that these functions have the 
Fourier expansions of the forms 
$f_r(z)=\sum_{m=1}^\infty c_m^{(r)}q^m$.

\begin{defn}[multiple Hecke $L$-function] \label{defL}
For $s\in \C$ with ${\mathrm{Re}\ s\gg 0}$ and for any
integers $\a_r\ge 1$, we define
\begin{align*}
L(s)&=L\left(
        \begin{smallmatrix}
        s,&\a_2,&\ldots,&\a_n\\
        f_1,&f_2,&\ldots,&f_n
        \end{smallmatrix}\right)
:=(-2\pi i)^{-(s+\a_2+\cdots+\a_n)}
\sum_{\substack{m_1>\cdots >m_n>0\\ (m_{n+1}:=0)}}
\frac{c_{m_1-m_2}^{(1)}\cdots
c_{m_{n}-m_{n+1}}^{(n)}}
{m_1^{s}m_2^{\a_2}\cdots m_n^{\a_n}}\\
&=(-2\pi i)^{-(s+\a_2+\cdots+\a_n)}
\sum_{l_1,\ldots,l_n>0}
\frac{c_{l_1}^{(1)}\cdots
c_{l_{n}}^{(n)}}
{(l_1+\cdots+l_n)^{s}(l_2+\cdots+l_n)^{\a_2}\cdots l_n^{\a_n}}.
\end{align*}
\end{defn}
Since Fourier coefficients $c_m^{(r)}$ have polynomial order in $m$ 
for any $r$,
there exists a constant $M>0$ such that
$\lvert c_{m_1-m_2}^{(1)}\cdots
c_{m_{n}-m_{n+1}}^{(n)} \rvert=O(m_1^M)$.
This implies the absolute convergence of $L(s)$ in $\mathrm{Re}(s)>M+1$.

In the following theorem,
we claim that $L(s)$ can be extended to an entire function on $\C$
for every fixed positive integers $\a_2,\ldots,\a_n$. However the analytic continuation of this kind of Dirichlet series had been 
studied
by Matsumoto-Tanigawa \cite{MT} in a general context.
They can regard other parameters $\a_2,\ldots,\a_n$ as complex
variables and show the analytic continuation to $\C^n$ by using
the Mellin-Barnes formula.
Our method in the next section is very simple and
different from theirs, but may not treat other parameters as variables.

\begin{thm}\label{th1}
For any fixed positive integers $\a_r$ for $2\le r \le n$,
the function L(s) can be extended holomorphically
to $\C$. Further, the following equations hold for any $s\in \C$: 
\begin{flalign}
&(i)\  I_{\ii}^{0}\left(\begin{smallmatrix}
s,&\a_2,&\ldots,&\a_n\\
f_1,&f_2,&\ldots,&f_n
\end{smallmatrix}\right)
=
\Gamma^{(s,\a_2,\ldots,\a_n)}\times \notag\\
&\qquad \qquad \qquad\sum_{\substack{0\le j_r< \a_r+j_{r+1}\\
(2\le \forall r \le n)\\ j_{n+1}:=0}}
\tbinom{s+j_2-1}{j_{2}}
\prod_{l=3}^{n}
\tbinom{\a_{l-1}+j_l-1}{j_{l}}
L\left(
\begin{smallmatrix}
  s+j_2,&\a_2-j_2+j_3,&\ldots,&\a_n-j_n+j_{n+1}&\\
  f_1,&f_2,&\ldots,&f_n
\end{smallmatrix}
\right),&\label{th11} \\
&(ii)\  L\left(\begin{smallmatrix}
s,&\a_2,&\ldots,&\a_n\\
f_1,&f_2,&\ldots,&f_n
\end{smallmatrix}\right)
=
\frac{1}{\Gamma^{(s,\a_2,\ldots,\a_n)}}\times \notag \\
&\qquad \qquad \qquad \quad\sum_{\substack{0\le j_r< \a_r\\
(2\le \forall r \le n)\\ j_{n+1}:=0}}
\tbinom{s+j_2-1}{j_{2}}
\prod_{l=2}^{n}
(-1)^{j_l}\tbinom{\a_{l}-1}{j_{l}}
I^0_{\ii}\left(
\begin{smallmatrix}
  s+j_2,&\a_2-j_2+j_3,&\ldots,&\a_n-j_n+j_{n+1}& \\
  f_1,&f_2,&\ldots,&f_n
\end{smallmatrix}
\right).&\label{th12}
 \end{flalign}
where $\Gamma^{(s,\a_2,\ldots,\a_n)}
=(-1)^{n}\Gamma(s)\Gamma(\a_2)\cdots\Gamma(\a_n)$.
\end{thm}
\begin{cor}\label{cor1}
For any integers $\a_r\ge 1$, we have
\begin{flalign*}
&(i)\ I_{\ii}^{0}\left(\begin{smallmatrix}
\a_1,&\ldots,&\a_n\\
f_1,&\ldots,&f_n
\end{smallmatrix}\right)
=
\Gamma^{(\a_1,\ldots,\a_n)}
\sum_{\substack{0\le j_r< \a_r+j_{r+1}\\
(2\le \forall r \le n)\\
\ j_1=j_{n+1}:=0}}
\prod_{l=2}^{n}
\tbinom{\a_{l-1}+j_l-1}{j_{l}}
L\left(
\begin{smallmatrix}
  \a_1-j_1+j_2,&\ldots,&\a_n-j_n+j_{n+1}\\
  f_1,&\ldots,&f_n
\end{smallmatrix}
\right),&\\
&(ii)\ L\left(
\begin{smallmatrix}
  \a_1,&\ldots,&\a_n\\
  f_1,&\ldots,&f_n
\end{smallmatrix}
\right)
=\frac{1}{\Gamma^{(\a_1,\ldots,\a_n)}}
\sum_{\substack{0\le j_r<\a_r\\
 (2\le \forall r \le n)\\ j_1=j_{n+1}:=0}}
\prod_{l=2}^{n} (-1)^{j_l}
\tbinom{\a_{l}-1}{j_{l}}
I_{\ii}^{0}\left(\begin{smallmatrix}
\a_1-j_1+j_2,&\ldots,&\a_n-j_n+j_{n+1}\\
f_1,&\ldots,&f_n
\end{smallmatrix}\right).&
\end{flalign*}
\end{cor}

\section{Proof}\label{sec3}

Corollary \ref{cor1} follows directly from Theorem \ref{th1}
by substituting $s=\a_1$. Hence we will prove Theorem \ref{th1}.
The procedure of the proof of (i) is as follows: 
The LHS of (i) is the Mellin
transformation of the function 
\begin{align*}
&F^{z}_{a} \left(
        \begin{smallmatrix}
       -, &\a_2,&\ldots,&\a_n\\
         f_1,&f_2,&\ldots,&f_n
        \end{smallmatrix}\right):=
f_1(z) I^{z}_{a} \left(
        \begin{smallmatrix}
        \a_2,&\ldots,&\a_n\\
        f_2,&\ldots,&f_n
        \end{smallmatrix}\right)
\end{align*}
for $a=\ii$.
In the first step, we define an alternative family of functions  $\widetilde{F}_{a}^z
(\begin{smallmatrix}
-, & \a_2,&\ldots,&\a_n\\
f_1,&f_2,&\ldots,&f_n
\end{smallmatrix})$ (Definition \ref{deftilF}) 
then describe $F_{a}^z$ by a linear combination 
of $\widetilde{F}_{a}^z$'s (Proposition \ref{FF}).
In the second step, we compute the Fourier expansion and Mellin transformation
of $\widetilde{F}_a^z$ (Proposition \ref{Fexp},
\ref{ME}). In particular, the Mellin transformation of 
$\widetilde{F}_a^z$ gives the analytic continuation of $L(s)$. 
Last, we will compute 
the Mellin transformation of ${F}_a^z$ by combining the results in 
the first and second steps. For the proof of (ii), 
we use an inversion description of $\widetilde{F}_{a}^z$ in terms of 
$F_{a}^z$'s in Proposition \ref{FF}.   

\begin{defn} \label{deftilF}
Let $\a_r$ be positive integers for $r=1,\ldots,n$.
For $a\in \overline{\H}$,
we define holomorphic functions
$\widetilde{I}^{z}_{a}$ and $\widetilde{F}^{z}_{a}$ in $z\in \overline{\H}$ by
\begin{align*}
\Iia:=&\int_{a}^{z}f_1(z_1)(z_1-z)^{\a_1-1}dz_1
\int_{a}^{z_1}
\cdots
\int_{a}^{z_{n-1}}f_n(z_n)(z_n-z_{n-1})^{\a_n-1}dz_n,
\end{align*}
\begin{align*}
\Fia:=
f_1(z)\widetilde{I}^{z}_{a} \left(
        \begin{smallmatrix}
        \a_2,&\ldots,&\a_n\\
        f_2,&\ldots,&f_n
        \end{smallmatrix}\right), \quad a\in \overline{\H}.
\end{align*}
\end{defn}
The integral converges again
by virtue of the cusp conditions of $f_r$'s.
In this definition, the parameters $\a_r$ should be positive integers,
otherwise the integral depends on the choice of path because the factor $(z_r-z_{r-1})^{\a_r-1}$ has a
singularity on $\H$ at $z_r=z_{r-1}$.
By definition $\widetilde{F}^{z}_{a}\binom{-}{f_1}=f_1(z)$
is independent on $a$.  The integral
$\widetilde{I}^{z}_{a}\tbinom{\a}{f}=
\int^{z}_{a}f(z_1)(z_1-z)^{\a-1}dz_1$
is known as the Eichler integral of $f$ when $f\in S_k(\G)$ and
$\a=k-1$. In this sense, $\widetilde{I}^{z}_{a}$ is an
iterated version of the Eichler integral. (A
modular property of $\widetilde{I}^{z}_{a}$ will be stated 
in the next section.)

\begin{prop}\label{FF} For $z,a\in \overline{\H}$, we have
\begin{flalign}
&(i)\ {F}^{z}_{a} \left(
\begin{smallmatrix}
  -,&\a_2,&\ldots,&\a_n\\
  f_1,&f_2&\ldots,&f_n
\end{smallmatrix}
\right)
=
\sum_{\substack{0\le j_r< \a_r+j_{r+1}\\
 (2\le \forall r\le n)\\
j_{n+1}:=0}}
\prod_{l=2}^{n}
\tbinom{\a_l+j_{l+1}-1}{j_l}
z^{j_2}
\widetilde{F}^{z}_{a} \left(
\begin{smallmatrix}
  -,&\a_2-j_2+j_3,&\ldots,&\a_n-j_n+j_{n+1}\\
  f_1,&f_2&\ldots,&f_n
\end{smallmatrix}
\right),&\label{FF1} \\ 
&(ii)\ {\widetilde{F}}^{z}_{a} \left(
\begin{smallmatrix}
  -,&\a_2,&\ldots,&\a_n\\
  f_1,&f_2&\ldots,&f_n
\end{smallmatrix}
\right)
=
\sum_{\substack{0\le j_r< \a_r\\
 (2\le \forall r\le n)\\
j_{n+1}:=0}}
\prod_{l=2}^{n}
(-1)^{j_l}\tbinom{\a_l-1}{j_l}
z^{j_2}
{F}^{z}_{a} \left(
\begin{smallmatrix}
  -,&\a_2-j_2+j_3,&\ldots,&\a_n-j_n+j_{n+1}\\
  f_1,&f_2&\ldots,&f_n
\end{smallmatrix}
\right).& \label{FF2}
\end{flalign}
\end{prop}
\pf Induction on $n$. For (i), when $n=2$ it is easy. For $n>2$, by using
the inductive hypothesis, we have
\begin{align*}
&{F}^{z}_{a} \left(
\begin{smallmatrix}
  -,&\a_2,&\ldots,&\a_n\\
  f_1,&f_2&\ldots,&f_n
\end{smallmatrix}
\right)
=f_1(z)\int^{z}_{a}z_2^{\a_2-1}
{F}^{z_2}_{a} \left(
\begin{smallmatrix}
  -,&\a_3,&\ldots,&\a_n\\
  f_2,&f_3&\ldots,&f_n
\end{smallmatrix}
\right)dz_2\\
&=
f_1(z)\int^{z}_{a}z_2^{\a_2-1}
\sum_{\substack{0\le j_r< \a_r+j_{r+1}\\
 (3\le \forall r\le n)\\
j_{n+1}:=0}}
\prod_{l=3}^{n}
\tbinom{\a_l+j_{l+1}-1}{j_l}
z_2^{j_3}
\widetilde{F}^{z_2}_{a} \left(
\begin{smallmatrix}
  -,&\a_3-j_3+j_4,&\ldots,&\a_n-j_n+j_{n+1}\\
  f_2,&f_2&\ldots,&f_n
\end{smallmatrix}
\right)dz_2.
\end{align*}
By the following binomial expansion
\begin{align*}
z_2^{\a_2+j_3-1}=(z_2-z+z)^{\a_2+j_3-1}=\sum_{0\le j_2<\a_2+j_3}
\tbinom{\a_2+j_{3}-1}{j_2}(z_2-z)^{\a_2-j_2+j_3-1}z^{j_2}
\end{align*}
we get (i). (ii) follows directly from the binomial expansions.  
\qed

\begin{prop}[Fourier expansion] \label{Fexp}
The Fourier expansions of $\widetilde{I}_{\ii}^z$ and $\widetilde{F}_{\ii}^z$ are given by
\begin{align*}
\Ii&=\frac{\Gamma^{(\a_1,\ldots,\a_n)}}
{(-2\pi i)^{\a_1+\cdots+\a_n}}\cdot
\sum_{\substack{m_1>\cdots >m_n>0\\
(m_{n+1}:=0)}}
\frac{c_{m_1-m_2}^{(1)}\cdots
c_{m_{n}-m_{n+1}}^{(n)}}
{m_1^{\a_1}\cdots m_n^{\a_n}}q^{m_1}, \\
\Fi&=\frac{\Gamma^{(\a_2,\ldots,\a_n)}}
{(-2\pi i)^{\a_2+\cdots+\a_n}}\cdot
\sum_{\substack{m_1>\cdots >m_n>0\\
(m_{n+1}:=0)}}
\frac{c_{m_1-m_2}^{(1)}\cdots
c_{m_{n}-m_{n+1}}^{(n)}}
{m_2^{\a_2}\cdots m_n^{\a_n}}q^{m_1}.
\end{align*}
In particular $\widetilde{I}_{\ii}^z$ and
$\widetilde{F}_{\ii}^z$
are periodic functions in $z$ of period one.
\end{prop}

\pf One can prove these by comparing the derivative
of both-hand sides in $z$, and by an inductive argument. \qed

\begin{prop}[Mellin transform]\label{ME}
\begin{align}
&\int_{\ii}^{0}\Ii z^{s-1}dz=
\Gamma^{(s,\a_1,\ldots,\a_n)}
L\left(
        \begin{smallmatrix}
        s+\a_1,&\a_2,&\ldots,&\a_n\\
        f_1,&f_2,&\ldots,&f_n
        \end{smallmatrix}\right),  \notag \\ 
&\int_{\ii}^{0}\Fi z^{s-1}dz=
\Gamma^{(s,\a_2,\ldots,\a_n)}
L\left(
        \begin{smallmatrix}
        s,&\a_2,&\ldots,&\a_n\\
        f_1,&f_2,&\ldots,&f_n
        \end{smallmatrix}\right).\label{ME2}
\end{align}
In particular, by eq. \eqref{ME2} we can extend $L(s)$ to an entire function.
\end{prop}
\pf By applying the Mellin transformation
to the equations in Proposition \ref{Fexp}, we easily get  
both expressions. Since the LHS of eq. \eqref{ME2} and 
$1/\Gamma(s)$ are entire, $L(s)$ can be extended entirely. 
\qed\\

\noindent{\bf Proof of Theorem \ref{th1}.} 
In Proposition \ref{ME}, we have already proved the 
analytic continuation of $L(s)$. For (i), by using eq. \eqref{FF1} and \eqref{ME2}, 
\begin{align*}
&I_{\ii}^{0}\left(\begin{smallmatrix}
s,&\a_2,&\ldots,&\a_n\\
f_1,&f_2,&\ldots,&f_n
\end{smallmatrix}\right)
=
\int_{\ii}^0 
F_{\ii}^{z}\left(\begin{smallmatrix}
-,&\a_2,&\ldots,&\a_n\\
f_1,&f_2,&\ldots,&f_n
\end{smallmatrix}\right)z^{s-1}dz\\
&=
\sum_{\substack{0\le j_r< \a_r+j_{r+1}\\
 (2\le \forall r\le n)\\
j_{n+1}:=0}}
\prod_{l=2}^{n}
\tbinom{\a_l+j_{l+1}-1}{j_l}
\int_{\ii}^0 
\widetilde{F}^{z}_{a} \left(
\begin{smallmatrix}
  -,&\a_2-j_2+j_3,&\ldots,&\a_n-j_n+j_{n+1}\\
  f_1,&f_2&\ldots,&f_n
\end{smallmatrix}
\right)z^{s+j_2-1}dz\\
&=
\sum_{\substack{0\le j_r< \a_r+j_{r+1}\\
 (2\le \forall r\le n)\\
j_{n+1}:=0}}
\prod_{l=2}^{n}
\tbinom{\a_l+j_{l+1}-1}{j_l}
\Gamma^{(s+j_2,\a_2-j_2+j_3,\ldots,\a_n-j_n+j_{n+1})}
L\left(
        \begin{smallmatrix}
        s+j_2,&\a_2-j_2+j_3,&\ldots,&\a_n-j_n+j_{n+1}\\
        f_1,&f_2,&\ldots,&f_n
        \end{smallmatrix}\right). 
\end{align*}
The elementary equation below completes the proof of (i): 
\begin{align*}
&\prod_{l=2}^{n}
\tbinom{\a_l+j_{l+1}-1}{j_l}
\Gamma^{(s+j_2,\a_2-j_2+j_3,\ldots,\a_n-j_n+j_{n+1})}
=
\Gamma^{(s,\a_2,\ldots,\a_n)}
\tbinom{s+j_{2}-1}{j_2}
\prod_{l=3}^{n}
\tbinom{\a_{l-1}+j_{l}-1}{j_l}. 
\end{align*}
For (ii), apply the Mellin transformation to 
eq. \eqref{FF2}, and use eq. \eqref{ME2}.
\qed


\section{Further properties of $\widetilde{I}_{a}^z$ and
$\widetilde{F}_{a}^z$}\label{sec4}
In this section, we show several further properties of
the functions $\widetilde{I}^{z}_{a}$ and $\widetilde{F}^{z}_{a}$ 
defined in the previous section.

\subsection{Alternative expression}

Since the function $\int_a^z f(w)(w-z)^{\a-1}dw$ in $z$ gives the $\a$-th anti-derivative of $f(z)$ up to a 
constant multiple in general, $\widetilde{I}_a^z$ 
has an alternative simple expression in terms of an iterated integral: 
\begin{prop} For any positive integers $\a_r$ we have 
\begin{align*}
&\Iia
=(-1)^{\a_1+\cdots+\a_n}\Gamma^{(\a_1,\ldots,\a_n)}\\
&\times \underbrace{\int_{a}^{z}dz_1\int_{a}^{z_1}\cdots
\int_{a}^{z_{\bu}}dz_{\bu}}_{\a_1-1}
\int_{a}^{z_\bu}f_1(z_\bu)dz_\bu
\cdots
\underbrace{\int_{a}^{z_\bu}dz_\bu\cdots
\int_{a}^{z_\bu}dz_{\bu}}_{\a_n-1}\int_{a}^{z_\bu}f_n(z_\bu)dz_\bu.
\end{align*}
\end{prop}

\subsection{Modular properties and its applications}

In this subsection we show modular properties of
$\widetilde{I}_{a}^z$ and $\widetilde{F}_{a}^z$
and discuss its applications.

Let $\a_r$ be positive integers for $r=1,\ldots,n$.
We fix the positive integers $k_r$ by
$k_r=\a_r+\a_{r+1}$ for $r=1,\ldots,n$ with $\a_{n+1}:=1$
throughout this subsection.
Conversely $\a_r$ are expressed by $k_r$'s as
$\a_r=(-1)^{n-r+1}+\sum_{j=r}^{n}(-1)^{j-r}k_j.$


\begin{prop} \label{modularity}
For $a\in \overline{\H}$ and $\g\in \mathrm{GL}^{+}_2(\Q)$
($+$ means positive determinant), we have
\begin{align*}
(I_n)\hspace{1cm}&\widetilde{I}^{z}_{a} \left(
        \begin{smallmatrix}
        \a_1,&\ldots,&\a_n\\
        f_1,&\ldots,&f_n
        \end{smallmatrix}\right)|_{(-\a_1+1)}\ \g =
\widetilde{I}^{z}_{\g^{-1}a} \left(
        \begin{smallmatrix}
        \a_1,&\ldots,&\a_n\\
        f_1|_{k_1}\g,&\ldots,&f_n|_{k_n}\g
        \end{smallmatrix}\right) \\
(F_n)\hspace{1cm}&\widetilde{F}^{z}_{a} \left(
        \begin{smallmatrix}
        -,&\a_2,&\ldots,&\a_n\\
        f_1,&f_2,&\ldots,&f_n
        \end{smallmatrix}\right)|_{(\a_1+1)} \g=
\widetilde{F}^{z}_{\g^{-1}a} \left(
        \begin{smallmatrix}
        -,&\a_2,&\ldots,&\a_n\\
        f_1|_{k_1} \g,&f_2|_{k_2} \g,&\ldots,&f_n|_{k_n} \g
        \end{smallmatrix}\right)
\end{align*}
where the slash action is defined in the usual manner:
$(f|_k \g)(z):=(\det \g)^{k/2}j(\g,z)^{-k}f(\g z)$
for any function $f$ on $\H$.
\end{prop}
\pf
First we check $(F_1)$ then show implications
$(F_n)\Rightarrow (I_n)$ and $(I_{n-1})\Rightarrow (F_n)$
for all $n$. $(F_1)$ is trivial:
$$f_1|_{(\a_1+1)}\g=f_1|_{k_1}\g$$
because of the choice of $k_1=\a_1+1$. For
$(F_n)\Rightarrow (I_n)$, we have
\begin{align}
\widetilde{I}^{z}_{a} &\left(
        \begin{smallmatrix}
        \a_1,&\ldots,&\a_n\\
        f_1,&\ldots,&f_n
        \end{smallmatrix}\right)|_{(-\a_1+1)}\ \g
=
\Big(\int_a^z
\widetilde{F}^{z}_{a} \left(
        \begin{smallmatrix}
        -,&\a_2,&\ldots,&\a_n\\
        f_1,&f_2,&\ldots,&f_n
        \end{smallmatrix}\right)(z_1-z)^{\a_1-1}dz_1\Big)
        \Big\arrowvert_{(-\a_1+1)}\g   \notag \\
&=(\det \g)^{\frac{-\a_1+1}{2}}j(\g,z)^{\a_1-1}
\int_a^{\g z}
\widetilde{F}^{z_1}_{a} \left(
        \begin{smallmatrix}
        -,&\a_2,&\ldots,&\a_n\\
        f_1,&f_2,&\ldots,&f_n
        \end{smallmatrix}\right)(z_1-\g z)^{\a_1-1}dz_1\notag \\
        &=(\det \g)^{\frac{-\a_1+1}{2}}j(\g,z)^{\a_1-1}
\int_{\g^{-1} a}^{z}
\widetilde{F}^{\g z_1}_{a} \left(
        \begin{smallmatrix}
        -,&\a_2,&\ldots,&\a_n\\
        f_1,&f_2,&\ldots,&f_n
        \end{smallmatrix}\right)(\g z_1-\g z)^{\a_1-1}d(\g z_1).
        \label{cont}
\end{align}
Substituting the formulas
\begin{align*}
\g z_1-\g z=(\det \g)j(\g,z_1)^{-1} j(\g,z)^{-1} (z_1-z), \quad
d(\g z_1)=(\det \g)j(\g,z_1)^{-2}dz_1
\end{align*}
and the induction hypothesis $(F_n)$:
\begin{align*}
\widetilde{F}^{\g z_1}_{a} \left(
        \begin{smallmatrix}
        -,&\a_2,&\ldots,&\a_n\\
        f_1,&f_2,&\ldots,&f_n
        \end{smallmatrix}\right)
=
(\det \g)^{-\frac{\a_1+1}{2}}j(\g,z)^{\a_1+1}
\widetilde{F}^{z_1}_{\g^{-1}a} \left(
        \begin{smallmatrix}
        -,&\a_2,&\ldots,&\a_n\\
        f_1|_{k_1} \g,&f_2|_{k_2} \g,&\ldots,&f_n|_{k_n} \g
        \end{smallmatrix}\right)
\end{align*}
into \eqref{cont}, we obtain the first implication:
\begin{align*}
\eqref{cont}
=
\int_{\g^{-1} a}^{z}
\widetilde{F}^{z_1}_{a} \left(
        \begin{smallmatrix}
        -,&\a_2,&\ldots,&\a_n\\
        f_1|_{k_1} \g,&f_2|_{k_2} \g,&\ldots,&f_n|_{k_n} \g
        \end{smallmatrix}\right) (z_1-z)^{\a_1-1}dz_1
=
\widetilde{I}^{z}_{\g^{-1}a} \left(
        \begin{smallmatrix}
        \a_1,&\ldots,&\a_n\\
        f_1|_{k_1}\g,&\ldots,&f_n|_{k_n}\g
        \end{smallmatrix}\right).
\end{align*}
For $(I_{n-1})\Rightarrow (F_n)$,
\begin{align}
\widetilde{F}^{z}_{a} \left(
        \begin{smallmatrix}
        -,&\a_2,&\ldots,&\a_n\\
        f_1,&f_2,&\ldots,&f_n
        \end{smallmatrix}\right)&|_{(\a_1+1)} \g
=
\Big(f_1(z)
\widetilde{I}^{z}_{a} \left(
        \begin{smallmatrix}
        \a_2,&\ldots,&\a_n\\
        f_2,&\ldots,&f_n
        \end{smallmatrix}\right)\Big)
\Big\arrowvert_{(\a_1+1)}\g \notag\\
&=(\det \g)^{\frac{\a_1+1}{2}}j(\g,z)^{-(\a_1+1)}
f_1(\g z)
\widetilde{I}^{\g z}_{a} \left(
        \begin{smallmatrix}
        \a_2,&\ldots,&\a_n\\
        f_2,&\ldots,&f_n
        \end{smallmatrix}\right) \label{cont2}.
\end{align}
By using the assumption $(I_{n-1})$:
\begin{align*}
\widetilde{I}^{\g z}_{a} \left(
        \begin{smallmatrix}
        \a_2,&\ldots,&\a_n\\
        f_2,&\ldots,&f_n
        \end{smallmatrix}\right)
=
(\det \g)^{-\frac{-\a_2+1}{2}}j(\g,z)^{-\a_2+1}
\widetilde{I}^{z}_{\g^{-1}a} \left(
        \begin{smallmatrix}
        \a_2,&\ldots,&\a_n\\
        f_2|_{k_2}\g,&\ldots,&f_n|_{k_n}\g
        \end{smallmatrix}\right),
\end{align*}
we have
\begin{align*}
\eqref{cont2}
&=
(\det \g)^{\frac{\a_1+\a_2}{2}}j(\g,z)^{-(\a_1+\a_2)}f_1(\g z)
\widetilde{I}^{z}_{\g^{-1}a} \left(
        \begin{smallmatrix}
        \a_2,&\ldots,&\a_n\\
        f_2|_{k_2}\g,&\ldots,&f_n|_{k_n}\g
        \end{smallmatrix}\right)\\
&=(f_1|_{(\a_1+\a_2)}\g)(z)
\widetilde{I}^{z}_{\g^{-1}a} \left(
        \begin{smallmatrix}
        \a_2,&\ldots,&\a_n\\
        f_2|_{k_2}\g,&\ldots,&f_n|_{k_n}\g
        \end{smallmatrix}\right)
=
\widetilde{F}^{z}_{\g^{-1}a} \left(
        \begin{smallmatrix}
        -,&\a_2,&\ldots,&\a_n\\
        f_1|_{k_1} \g,&f_2|_{k_2} \g,&\ldots,&f_n|_{k_n} \g
        \end{smallmatrix}\right).
\end{align*}
This completes the proof.
\qed \\
We show two applications of Proposition \ref{modularity}.

\subsubsection{Functional equations}
Let $w_N=\sm 0 & -1 \\ N & 0\esm$ (Fricke involution).
Suppose that $f_r$ are $w_N$-eigenfunctions of weight
$k_r$:
$f_r|_{k_r} w_N=\e_r f_r$,
$\e_r\in \{\pm 1\}$ for $r=1,\ldots,n$.
Put
$g(z):=\Iic$ for $c=\frac{i}{\sqrt{N}}$
($w_N$-fixed point) and set
$
\Lambda(s, g):=\int_{\ii}^{0}g(z)  z^{s-1}dz.
$

\begin{thm} It holds that
$$\Lambda(s,g)= (-1)^{\a_1-1+s}
N^{\frac{-\a_1+1}{2}-s}\e \Lambda(-\alpha_1+1-s,g)$$
where $\e=\prod_{r=1}^n\e_r$.
\end{thm}

{\pf}
For any holomorphic function $f(z)$ on $\H$ satisfying the
condition (ii) in \S 1, one can show by the standard method
that the Mellin transformation
$\Lambda(s, f):=\int_{\ii}^{0}f(z)  z^{s-1}dz$ satisfies
\begin{align*}
\Lambda(s,\widehat{f})=(-1)^{s-\a}N^{\frac{\a}{2}-s}
\Lambda(\a-s,f)
\end{align*}
where $\a$ is any integer and $\widehat{f}(z)=(f|_{\a}w_N)(z)$.
If we substitute $g$ to $f$ and $-\a_1+1$ to $\a$, then we obtain the theorem because $\widehat{g}=\e g$ is true by virtue of Proposition
\ref{modularity}.
\qed\\


\subsubsection{Twisted iterated integral}

For any Dirichlet character $\chi$ mod $M$ we set
$f^{\chi}(z):=\sum_{m=1}^{M}\chi(m)f(\frac{z+m}{M}).$
It is well-known \cite{R} that a holomorphic function $f$ on
$\mathbf{H}$ is $\Gamma_0(N)$-invariant of weight $k$:
$f|_k \g=f$ for $\forall \g \in \Gamma_0(N)$
if and only if it is
periodic with period one and satisfies
$$
f^{\chi}\big|_k\s =\chi(-1)
f^{\overline{\chi}},
\qquad \s=\sm 0 & -1 \\ 1 & 0\esm \label{twist}
$$
for all Dirichlet characters $\chi
\pmod{Nc}$ for each $c\in \{1,\ldots,N\}$, where
$\overline{\chi}$ is the complex conjugate of $\chi$.

Suppose $f_r\in S_{k_r}(\G_0(N))$ and $\chi_r$ are Dirichlet
characters of mod $Nc$ with $1\le c\le N$
for $r=1,\ldots,n$.
Define two functions by
\begin{align*}
G_{\chi}(z):&=\Iichi, \qquad
G_{\overline{\chi}}(z):=\Iichibar
\end{align*}
and denote their Mellin transformations by
$\Lambda(s,G_\chi)$ and $\Lambda(s,G_{\overline{\chi}})$
respectively.
By Proposition \ref{modularity} and above fact, one can show that
\begin{align*}
G_{\chi}\big|_{(-\a_1+1)}\s=\chi_1(-1)\cdots \chi_n(-1)
G_{\overline{\chi}}.
\end{align*}
Hence we have the following identity:
\begin{thm}
\begin{align*}
\Lambda(s,G_\chi)
=(-1)^{s} \chi_1(-1)\cdots \chi_n(-1)
\Lambda(1-\a_1-s,G_{\overline{\chi}}).
\end{align*}
\end{thm}



\end{document}